\documentclass[10pt, a4paper, reqno]{amsart}

\usepackage{amsmath,amssymb,amsthm}
\usepackage{enumitem}
\usepackage{xcolor}
\usepackage[colorlinks=true,linkcolor=blue,citecolor=blue,urlcolor=blue]{hyperref}

\numberwithin{equation}{section}

\newtheorem{theorem}{Theorem}[section]
\newtheorem{lemma}[theorem]{Lemma}
\newtheorem{proposition}[theorem]{Proposition}
\newtheorem{corollary}[theorem]{Corollary}
\newtheorem*{thmA}{Theorem A}
\newtheorem*{thmB}{Theorem B}
\newtheorem*{thmC}{Theorem C}

\theoremstyle{definition}
\newtheorem{definition}[theorem]{Definition}

\theoremstyle{remark}
\newtheorem{remark}[theorem]{Remark}

\newcommand{\CC}{\mathbb{C}}
\newcommand{\RR}{\mathbb{R}}
\newcommand{\NN}{\mathbb{N}}
\newcommand{\ZZ}{\mathbb{Z}}
\newcommand{\DD}{\mathbb{D}}
\newcommand{\TT}{\mathbb{T}}

\newcommand{\HH}{\mathbf{H}_{q}}
\newcommand{\HHk}[1]{\mathbf{H}_{q}^{(#1)}}
\newcommand{\Fock}{\mathcal{F}}
\newcommand{\Lq}{L^{2}(\CC,\mu_{q})}
\newcommand{\Poly}{\CC[z,\zb]}
\newcommand{\Pgr}{\mathcal{P}}

\newcommand{\Rq}{R_{q}}
\newcommand{\Rqb}{\overline{R}_{q}}
\newcommand{\Mz}{M_{z}}
\newcommand{\Mzb}{M_{\zb}}
\newcommand{\Lam}{\Lambda_{q}}
\newcommand{\Nz}{\widehat{N}_{z}}
\newcommand{\Nzb}{\widehat{N}_{\zb}}
\newcommand{\J}{\widehat{J}}

\newcommand{\muq}{d\mu_{q}}
\newcommand{\Eq}{E_{q}}
\newcommand{\Aq}{A_{q}}
\newcommand{\Bq}{B_{q}}
\newcommand{\zb}{\overline{z}}
\newcommand{\wb}{\overline{w}}
\newcommand{\pp}{\partial}
\newcommand{\ppb}{\overline{\partial}}
\newcommand{\qbin}[2]{\genfrac{[}{]}{0pt}{}{#1}{#2}_{q}}
\newcommand{\qint}[1]{[#1]_{q}}
\newcommand{\qfac}[1]{[#1]_{q}!}
\newcommand{\qiint}[1]{[#1]_{q^{-1}}}
\newcommand{\rad}{\rho_{q}}

\title[Polyanalytic $q$-Fock spaces]
  {Polyanalytic $q$-Fock spaces and the polyanalytic $q$-Bargmann transform}

\author[I. Ali]{Irfan Ali}
\address{
 Irfan Ali:
  \endgraf
Department of Mathematics: Analysis, Logic and Discrete Mathematics
  \endgraf
 Ghent University
 Belgium
   \endgraf
  and
  \endgraf
  Department of Mathematical Science,
  BUITEMS, Pakistan}
\email{irfan.ali@ugent.be}

\author[U.\,K\"{a}hler]{Uwe K\"{a}hler}
\address{
  Uwe K\"{a}hler:
  \endgraf
  CIDMA, Departamento de Matem\'{a}tica, Universidade de Aveiro,
  \endgraf
  Campus de Santiago, P-3810-193, Aveiro,
  Portugal}
\email{ukaehler@ua.pt}

\author[M. Ruzhansky]{Michael Ruzhansky}
\address{
 Michael Ruzhansky:
  \endgraf
Department of Mathematics: Analysis, Logic and Discrete Mathematics
  \endgraf
 Ghent University
 Belgium
 \endgraf
   and
  \endgraf
  School of Mathematical Sciences, Queen Mary University of London
  \endgraf
  London, United Kingdom.}
\email{michael.ruzhansky@ugent.be}



\begin{document}

\begin{abstract}

We study polyanalytic $q$-Fock spaces from an operator-theoretic point of view. We compute the formal adjoint $A^*$ of the Jackson derivative $A=\Rqb$ on the polynomials in $z,\zb$ inside $\Lq$ and we show that, up to unimodular factors, a twist by the grading operator is the only modification of $A$ diagonal in the grading for which a deformed canonical relation $aa^*-\lambda a^*a=I$ holds. The resulting second-order $q$-difference operator $a^*a$ is essentially self-adjoint on the polynomials and $q$-complex Hermite polynomials span its eigenspaces. These eigenspaces give the orthogonal decomposition of $\Lq$ into true polyanalytic $q$-Fock spaces. We determine the disk on which each of them is a reproducing kernel Hilbert space and we factor the corresponding unitary $q$-Bargmann transforms through powers of $a^*$.
\end{abstract}

\maketitle

\section{Introduction}\label{sec:intro}

Classical true polyanalytic Fock spaces can be characterized as eigenspaces of the Landau Hamiltonian (magnetic Laplacian); see, e.g.,~\cite{Mouayn, RT16, RT17a, RT17b, RT19a, RT19b, Vas}. In fact, the orthogonal decomposition of the polyanalytic Fock space corresponds to the decomposition into pure Landau level eigenspaces. Their applications are not restricted to physics only but also include time-frequency analysis and in particular, sampling and interpolation results for vector-valued Gabor systems with Hermite windows~\cite{AbreuACHA,AF}. The Bargmann transform~\cite{Barg}, which maps $L^{2}(\RR)$ unitarily onto the Fock space $\Fock^{2}$ of entire functions square integrable against $\pi^{-1}e^{-|z|^{2}}dxdy$, carrying the Hermite basis $\psi_{m}$ to $\frac{z^{m}}{\sqrt{m!}}$. In $L^{2}(\CC,\pi^{-1}e^{-|z|^{2}}dxdy)$ the operators $A=\ppb$ and $A^*=\Mzb-\pp$ are mutually adjoint on polynomials and satisfy $[A,A^*]=I$. For each $n\in\NN_{0}$, the true polyanalytic Fock space $\Fock^{(n)}$ coincides with the eigenspace of the operator $A^{*}A$ associated with the eigenvalue $n$. Moreover, each $\Fock^{(n)}$ admits an orthonormal basis given by the normalized complex Hermite polynomials introduced by It\^{o}~\cite{Ito}.

The work of Z.~Mouayn et~al.~\cite{Mouayn} primarily focuses on the explicit construction of reproducing kernels and coherent states. In connection with this $q$-deformed oscillator, Arik and Coon~\cite{Arik} introduced a one-parameter family of reproducing kernel Hilbert spaces $\HH$ of analytic functions, with kernel $\Eq(z\wb)$, governed by the relation $\Rq\Mz-q\Mz\Rq=I$, where $\Rq$ is the Jackson derivative and $\Mz$ is multiplication by $z$. These spaces are uniquely determined by the requirement $\Mz^*=\Rq$ and they interpolate between the Hardy space of the disk obtained by setting $q=0$ in the defining formulas and the Fock space obtained formally as $q\to1^{-}$~\cite{ACKS,AlpayPorat}. They are realised inside $\Lq$ for a measure $\mu_{q}$ which is discrete in the radial direction and it is carried by countably many circles $|z|=r_{j}$, $r_{j}=q^{\frac{j}{2}}(1-q)^{-\frac{1}{2}}$.

Our main results on polyanalytic $q$-Fock spaces are as follows. We compute the adjoint of $A=\Rqb$ on $\Poly\subset\Lq$ and show that, up to unimodular factors, the twist by the grading operator is the only modification of $A$, diagonal in the grading by $\deg_{\zb}-\deg_{z}$, for which a deformed canonical relation $aa^*-\lambda a^*a=I$ holds; necessarily $\lambda=q^{-1}$ (Theorem~A). The polynomials $H_{m,n}=q^{mn}(A^*)^{n}z^{m}$ generated by the formal adjoint are $q$-analogues of It\^{o}'s complex Hermite polynomials $(\Mzb-\pp)^{n}z^{m}$~\cite{Ito}: $A$ and $A^*$ act on them as lowering and raising operators, they form an orthogonal basis of $\Poly$ consisting of eigenvectors of $a^*a$ and their expansion and generating function reduce formally to It\^{o}'s as $q\to1^{-}$ (Theorem~B). Further, $\Lam=a^*a$, the analogue of $A^*A$, is essentially self-adjoint on $\Poly$ with spectrum $\{\qiint{n}:n\in\NN_{0}\}$, each eigenvalue being independent of the analytic degree and of infinite multiplicity (Corollary~\ref{cor:esa}); its eigenspaces $\HHk{n}$ give the orthogonal decomposition of $\Lq$ into true $q$-polyanalytic Fock spaces (Theorem~C(i)), a $q$-analogue of Vasilevski's decomposition of $L^{2}(\CC,\pi^{-1}e^{-|z|^{2}}dxdy)$ into the spaces $\Fock^{(n)}$~\cite{Vas}. $\HHk{n}$ is a reproducing kernel Hilbert space of functions polyanalytic of order $n+1$ on $\DD_{r_n}$ and on no larger disk, with the norm of $\Lq$, which is not equivalent to that of $L^{2}(\DD_{r_n},\mu_{q})$ (Theorems~C(ii) and~D). All levels are generated from the analytic one by the single operator $a^*$: the unitary $q$-Bargmann transform onto $\HHk{n}$ is $(a^*)^{n}\Bq/\sqrt{\qiint{n}!}$ (Theorem~C(iii)) and the twist of Theorem~A is what makes one normalising constant serve for the whole basis. The $n+1$ circles $|z|=r_{0},\dots,r_{n}$ outside $\DD_{r_n}$ or on its boundary form a sampling family for $\HHk{n}$, whereas the circles inside do not (Theorem~D); in $\Fock^{(n)}$ no finite family of circles has this property.

\section{Preliminaries}\label{sec:prelim}

Throughout the paper $q\in(0,1)$ and $\NN_{0}=\{0,1,2,\dots\}$. We denote
\begin{equation*}
  \qint{k}=\frac{1-q^{k}}{1-q}\ \ (k\in\ZZ),\qquad
  \qfac{k}=\qint{1}\cdots\qint{k}=\frac{(q;q)_{k}}{(1-q)^{k}}\ \ (k\in\NN_{0}),
\end{equation*}
where $(a;q)_{k}=\prod_{0\le j<k}(1-aq^{\,j})$ for $k\in\NN_{0} {}\cup\{\infty\}$ and we recall that
\begin{equation}\label{eq:reciprocal}
  \qiint{k}=q^{\,1-k}\qint{k},
  \quad
  \qiint{k+1}-\qiint{k}=q^{-k},
  \quad
  \qiint{k}\nearrow\infty .
\end{equation}
Furthermore, we have
\begin{equation}\label{eq:qbinfac}
  \qbin{k}{m}=\frac{\qfac{k}}{\qfac{m}\qfac{k-m}}
  \quad(0\le m\le k).
\end{equation}
The $q$-analogues of the exponential function are given by
\begin{align}
\label{eq:qexp}
  \Eq(x)&=\sum_{k\ge0}\frac{x^{k}}{\qfac{k}}
   =\prod_{j\ge0}\frac{1}{1-x(1-q)q^{\,j}},
   \qquad |x|<\frac{1}{1-q},\\
  \Eq(x)^{-1}&=\sum_{k\ge0}
   \frac{(-1)^{k}q^{\binom{k}{2}}x^{k}}{\qfac{k}}
   =\bigl((1-q)x;q\bigr)_{\infty}\notag ,
\end{align}
where the second series converges for every $x\in\CC$. The $q$-analogue of Gauss' formula, for commuting $x$ and $a$, is given by
\begin{equation}
\label{eq:gauss}
  \prod_{l=0}^{k-1}\bigl(x+a q^{\,l}\bigr)
  =\sum_{i=0}^{k}\qbin{k}{i}q^{\binom{i}{2}}a^{i}x^{\,k-i}.
\end{equation}
The $q$-Jackson derivatives are defined by (see \cite{jackson})
\[
  (\Rq f)(z)=\frac{f(z)-f(qz)}{(1-q)z},
\]
and we use it, together with its conjugate $\Rqb$, as a replacement for the classical Wirtinger derivatives. Like in the classical case, $\Rq$ and $\Rqb$ are acting on $z$ and $\zb$ as independent variables. Explicitly, a function is written as $f=f(z,\zb)$ in two independent variables and
\begin{equation}
\label{eq:jackson2var}
  (\Rq f)(z,\zb)=\frac{f(z,\zb)-f(qz,\zb)}{(1-q)z},\quad
  (\Rqb f)(z,\zb)=\frac{f(z,\zb)-f(z,q\zb)}{(1-q)\zb},
\end{equation}
 so that on monomials $\Rq z^{a}\zb^{\,b}=\qint{a}z^{a-1}\zb^{\,b}$ and $\Rqb z^{a}\zb^{\,b}=\qint{b}z^{a}\zb^{\,b-1}$.

We denote $e_{a,b}=z^{a}\zb^{\,b}$ for $a,b\in\NN_{0}$  with $e_{a,b}=0$ if $a<0$ or $b<0$ and let $\Mz,\Mzb$ denote multiplication by $z,\zb$. 

The Jackson integral is defined as (see \cite{jackson})
\[
  \int_{0}^{c}h(s)\,d_{q}s=c(1-q)\sum_{j\ge0}q^{\,j}h(cq^{\,j}).
\]

Moreover, we denote
\begin{equation}
\label{eq:radii}
  \rad=(1-q)^{-\frac{1}{2}},\qquad
  r_{j}=q^{\frac{j}{2}}\rad,\qquad
  s_{j}=r_{j}^{2}=\frac{q^{\,j}}{1-q},\qquad j\in\NN_{0},
\end{equation}
and let $d\lambda_{r}$ be the normalised Lebesgue measure on $\TT_{r}:=\{z\in\mathbb{C}:|z|=r\}$ we also write $\DD_{r}=\{z\in\CC:|z|<r\}$, $\overline{\DD}_{r}$ for its closure. We then consider the measure (see \cite{LM})
\begin{equation}
\label{eq:measure}
  \muq=(q;q)_{\infty}\sum_{j\ge0}\frac{q^{\,j}}{(q;q)_{j}}\,
   d\lambda_{r_{j}} ,
\end{equation}
that is, $\mu_{q}=\sum_{j\ge0}w_{j}\lambda_{r_{j}}$ with $w_{j}=\frac{q^{\,j}(q;q)_{\infty}}{(q;q)_{j}}>0$. Since $r_{j}\downarrow0$, the support of $\mu_{q}$ is the compact set $K_{q}=  \{0\}\cup\bigcup_{j\ge0}\TT_{r_{j}}\subset \overline{\DD}_{\rad}$. In particular, every polynomial in $z,\zb$ belongs to $\Lq$, which carries the inner product $\langle f,g\rangle=\int_{\CC}f\,\overline{g}\,\muq$, linear in the first argument. Since $\mu_{q}$ restricted to $\TT_{r_j}$ is $w_{j}\lambda_{r_j}$ with $w_{j}>0$, the restriction of $f\in\Lq$ to $\TT_{r_j}$ is a well-defined element of $L^{2}(\TT_{r_j},d\lambda_{r_j})$; by monotone convergence $\|f\|^{2}=\sum_{j}w_{j}\|f|_{\TT_{r_j}}\|^{2}_{L^{2}(d\lambda_{r_j})}$, so that $\|f|_{\TT_{r_j}}\|_{L^{2}(d\lambda_{r_j})}\le w_{j}^{-1/2}\|f\|$. Moreover, we have the following lemma for the measure defined in \eqref{eq:measure}.

\begin{lemma}[cf.~\cite{LM}]\label{lem:moment}
We have $\int_{\CC}z^{a}\zb^{\,b}\muq=\qfac{a}\,\delta_{a,b}$ for all $a,b\in\NN_{0}$. Consequently, $\mu_{q}$ is a probability measure and
\begin{equation}\label{eq:gram}
  \langle e_{a,b},e_{c,d}\rangle=\int_{\CC}z^{a+d} \, \zb^{\,b+c}\muq=\delta_{a-b,\,c-d}\,\qfac{a+d},
  \qquad a,b,c,d\in\NN_{0}.
\end{equation}
\end{lemma} 

\begin{proof} 
By using \eqref{eq:measure}, the Euler's identity $\sum_{j\ge0}\frac{z^{\,j}}{(q;q)_{j}}=\frac{1}{(z;q)_{\infty}}$ and $(q^{\,a+1};q)_{\infty}=\frac{(q;q)_{\infty}}{(q;q)_{a}}$ we obtain
\begin{align*}
  \int_{\CC}z^{a}\zb^{\,b}\muq
  &=\delta_{a,b}\frac{(q;q)_{\infty}}{(1-q)^{a}}
   \sum_{j\ge0}\frac{\bigl(q^{\,a+1}\bigr)^{j}}{(q;q)_{j}}
  =\delta_{a,b}\frac{(q;q)_{\infty}}{(1-q)^{a}}
   \cdot\frac{1}{(q^{\,a+1};q)_{\infty}} =\delta_{a,b}\qfac{a}.
\end{align*}
The case $a=b=0$ shows that $\mu_{q}$ is a probability measure and identity \eqref{eq:gram} follows from $\langle e_{a,b},e_{c,d}\rangle=\int z^{a+d}\zb^{\,b+c}\muq$ together with the equivalence of $a+d=b+c$ and $a-b=c-d$.
\end{proof}

\begin{lemma}\label{lem:inj}
The map $\Poly\to\Lq$ is injective: if $f\in\Poly$ vanishes $\mu_{q}$-almost everywhere, then $f=0$. In particular, the monomials $e_{a,b}$ are linearly independent in $\Lq$.
\end{lemma}
\begin{proof}
Write $f=\sum c_{a,b}e_{a,b}$. Since $e_{a,b}(re^{i\theta})=r^{a+b}e^{-i(b-a)\theta}$,
\[
  f(re^{i\theta})=\sum_{J\in\ZZ}e^{-iJ\theta}\sum_{b-a=J}c_{a,b}\,r^{a+b}.
\]
As $w_{j}>0$, $f$ vanishes $\lambda_{r_j}$-almost everywhere on $\TT_{r_j}$. For each $J$ the polynomial $r\mapsto\sum_{b-a=J}c_{a,b}r^{a+b}$ vanishes at the infinitely many points $r_{j}$, so it is zero. Since $a+b=2a+J$ determines $(a,b)$ for fixed $J$, all $c_{a,b}$ vanish.
\end{proof}

Now, we introduce some definitions for our analysis of $q$-polyanalytic functions.

\begin{definition}\label{def:grading}
Let $\Nz$ and $\Nzb$ be the number operators on $\Poly$, $\Nz e_{a,b}=a\,e_{a,b}$ and $\Nzb e_{a,b}=b\,e_{a,b}$ and let $\J:=\Nzb-\Nz$ be the \emph{ grading operator}, so that $\J e_{a,b}=(b-a)\,e_{a,b}$. For $J\in\ZZ$ let
\[
  \Pgr_{J}=\operatorname{span}\{e_{a,b}:b-a=J\}
\]
be the \emph{$J$ sector}, i.e.\ the eigenspace of $\J$ with eigenvalue $J$. Then
\[
  \Poly=\bigoplus_{J\in\ZZ}\Pgr_{J}
\]
(algebraic direct sum) is the grading of $\Poly$ by $\deg_{\zb}-\deg_{z}$. For a function $\varphi:\ZZ\to\CC$ we define $\varphi(\J)$ on $\Poly$ by
\[
  \varphi(\J)\,e_{a,b}=\varphi(b-a)\,e_{a,b},
\]
so that $\varphi(\J)$ acts on $\Pgr_{J}$ as multiplication by $\varphi(J)$. In particular, for $s\in\RR$,
\[
  q^{s\J}e_{a,b}=q^{s(b-a)}e_{a,b}.
\]
\end{definition}
By \eqref{eq:gram} the components $\Pgr_{J}$ are mutually orthogonal in $\Lq$, whereas inside $\Pgr_{J}$, writing $b=a+J$ and $d=c+J$ with $a,c\ge\max(0,-J)$, one has $\langle e_{a,b},e_{c,d}\rangle =\qfac{a+c+J}$, so that the Gram matrix is the Hankel matrix $\bigl(\qfac{a+c+J}\bigr)_{a,c\ge\max(0,-J)}$.  Since the operator $\J$ is diagonal with real eigenvalues and the $\Pgr_{J}$ are orthogonal, every real-valued function of $\J$ is symmetric on $\Poly$. This is used for $q^{\pm\frac{\J}{2}}$ without further comment.

The classical $q$-Fock space, or Arik--Coon space, $\HH$ is the closure of $\CC[z]$ in $\Lq$. By Lemma~\ref{lem:moment} the monomials are orthogonal in $\HH$ with $\|z^{m}\|^{2}=\qfac{m}$, so $e_{m,q}:=\frac{z^{m}}{\sqrt{\qfac{m}}}$ forms an orthonormal basis of $\HH$.  Consequently,  every $f\in\HH$ satisfies $f=\sum_{m}\langle f,e_{m,q}\rangle e_{m,q}$ in $\Lq$ and writing $\langle f,e_{m,q}\rangle=a_{m}\sqrt{\qfac{m}}$, The Parseval identity gives
\[
  \|f\|^{2}=\sum_{m}|a_{m}|^{2}\qfac{m}.
\]
 By using the Cauchy--Schwarz inequality, we obtain
\[
  \sum_{m}|a_{m}|\,|z|^{m}
  \;\le\;
  \Bigl(\sum_{m}|a_{m}|^{2}\qfac{m}\Bigr)^{\frac12}
  \Bigl(\sum_{m}\frac{|z|^{2m}}{\qfac{m}}\Bigr)^{\frac12}
  \;=\;
  \|f\|\,\Eq\bigl(|z|^{2}\bigr)^{\frac12},
\]
for every $z$ at which $\Eq(|z|^{2})$ is finite for $|z|<\rad$. 

$\HH$ is a reproducing kernel Hilbert space of analytic functions on $\DD_{\rad}$ with kernel $\Eq(z\wb)$. Moreover, the inner product of $\HH$ is the only inner product on $\CC[z]$ with $\|1\|=1$ for which $\Mz$ and $\Rq$ are formal adjoints of each other~\cite{ACKS, Arik}: such an inner product satisfies $\langle z^{m},z^{n}\rangle=\langle z^{m-1},\Rq z^{n}\rangle=\qint{n}\langle z^{m-1},z^{n-1}\rangle$ for $m\ge1$ and $\langle1,z^{n}\rangle=\overline{\langle z^{n-1},\Rq1\rangle}=0$ for $n\ge1$, hence $\langle z^{m},z^{n}\rangle=\delta_{m,n}\qfac{n}$. Setting $q=0$ these gives the classical Hardy space $H^{2}(\DD)$ and formally, as $q\to1^{-}$, one obtains the classical Fock space $\Fock^{2}$; in this sense the spaces $\HH$ interpolate between the Hardy and the Fock space.

We will also need the following lemma.

\begin{lemma}\label{lem:polar}
For $f,g\in\Lq$ we have
\begin{equation}
\label{eq:polar}
  \langle f,g\rangle=\frac{1}{2\pi}\int_{0}^{2\pi}\!\!
   \int_{0}^{\frac{1}{1-q}}
   f\bigl(\sqrt{s}e^{i\theta}\bigr)
   \overline{g\bigl(\sqrt{s}e^{i\theta}\bigr)}\,
   \frac{d_{q}s}{\Eq(qs)}\,d\theta .
\end{equation}
In particular $\sum_{j\ge0}w_{j}s_{j}^{\,p}=\qfac{p}$ for every $p\in\NN_{0}$.
\end{lemma}

\begin{proof}
By the definition of the Jackson integral,
\[
\int_0^{\frac{1}{1-q}} h(s)\,d_qs
=
\sum_{j=0}^{\infty}
q^j
h\left(\frac{q^j}{1-q}\right)
=
\sum_{j=0}^{\infty}
q^j
h(s_j), \,\, \text{where} \,\, s_j=\frac{q^j}{1-q}.
\]
Moreover, by \eqref{eq:qexp} we get
\[
\frac{1}{\Eq(qs_j)}
=
\bigl((1-q)qs_j;q\bigr)_\infty
=
(q^{j+1};q)_\infty
=
\frac{(q;q)_\infty}{(q;q)_j}.
\]
We can write form \eqref{eq:polar} as
\[
\begin{aligned}
&\frac{1}{2\pi}
\int_0^{2\pi}
\int_0^{\frac{1}{1-q}}
f\bigl(\sqrt{s}e^{i\theta}\bigr)
\overline{g\bigl(\sqrt{s}e^{i\theta}\bigr)}
\frac{d_qs}{\Eq(qs)}\,d\theta
\\
&\qquad=
\sum_{j=0}^{\infty}
\frac{q^j(q;q)_\infty}{(q;q)_j}
\frac{1}{2\pi}
\int_0^{2\pi}
f\bigl(\sqrt{s_j}e^{i\theta}\bigr)
\overline{g\bigl(\sqrt{s_j}e^{i\theta}\bigr)}
\,d\theta
\\
&\qquad=
\sum_{j=0}^{\infty}
w_j
\int_{\mathbb C}
f(z)\overline{g(z)}\,d\lambda_{r_j}(z)
=
\int_{\mathbb C}
f(z)\overline{g(z)}\,d\mu_q(z)
=
\langle f,g\rangle.
\end{aligned}
\]
Here $f(\sqrt{s}e^{i\theta})$ is used only at $s=s_{j}$, where it is the restriction of $f$ to $\TT_{r_j}$. Taking $f=g=z^{p}$ and using Lemma~\ref{lem:moment} gives$\sum_{j\ge0}w_{j}s_{j}^{\,p}=\qfac{p}$, which is Jackson's $q$-Gamma integral $\Gamma_{q}(p+1)=\qfac{p}$.
\end{proof}
We can now define the $q$-polyanalytic functions as follows.
\begin{definition}\label{def:qpoly}  
Let $n\in\NN$ and let $\Omega\subset\CC$. A $q$-polyanalytic function of order $n$ is a function $f$ on a domain in the complex plane such that
\[
  \Rqb^{\,n}f=0 .
\]
\end{definition}

We are going to use the fact that a polynomial in $\zb$ with analytic coefficients satisfies $\Rqb^{\,n}f=0$ on a domain $\Omega$ if and only if
\begin{equation}\label{eq:almansi}
  f=\sum_{\nu=0}^{n-1}\zb^{\,\nu}f_{\nu},\qquad
  f_{\nu}\ \text{analytic on }\Omega ,
\end{equation}
the $f_{\nu}$ being unique. Indeed $\Rqb^{\,n}f=\sum_{\nu\ge n}\qint{\nu}\qint{\nu-1}\cdots\qint{\nu-n+1} \zb^{\,\nu-n}f_{\nu}$ and $\qint{j}\neq0$ for $j\ge1$. In the classic case, where $\Rqb^{\,n}f=0$ is replaced by $\ppb^{\,n}f=0$, formula \eqref{eq:almansi} is also known as the Almansi decomposition.

    In the next section, we present our main theorems. 
\section{Main Results}\label{sec:main}

The first theorem provides a characterization of the adjoint of the operator $A=\Rqb$ with respect to the inner product \eqref{eq:polar} and it also identifies the associated second-order operator.

\begin{thmA}
Let $A=\Rqb$ on $\Poly\subset\Lq$.
\begin{enumerate}[label=\textup{(\roman*)},leftmargin=*]
\item The adjoint of the operator $A$ on $\Poly$ is given by $A^*=(\Mzb-\Rq)\,q^{-\Nz}$. Moreover, the commutator satisfies $[A,A^*]=q^{\J}$, where $\J=\Nzb-\Nz$ denotes the grading operator given in Definition~\ref{def:grading}.
\item Let $D=d(\J)$ be diagonal in the grading with $d(J)\neq0$ for all $J$ and put $\widetilde a=DA$ and $\widetilde a^*=A^*\overline{d(\J)}$. If $\widetilde a\widetilde a^* -\lambda\widetilde a^*\widetilde a=I$ on $\Poly$ for some $\lambda\in\CC$, then $\lambda=q^{-1}$ and $|d(J)|^{2}=q^{-J}$ for all $J$.
Conversely, $a:=q^{-\frac{\J}{2}}A$ and $a^*:=A^*q^{-\frac{\J}{2}}$ are mutually adjoint on $\Poly$ and satisfy the following relation
\begin{equation}
\label{eq:AC}
  a\,a^*-q^{-1}a^*a=I .
\end{equation}
\item The operator $\Lam:=a^*a$ is symmetric on $\Poly$ and acts by $\Lam e_{a,b}=\qiint{b}\bigl(e_{a,b}-\qint{a}e_{a-1,b-1}\bigr).$ Its point spectrum on $\Poly$ is $$\sigma_{p}(\Lam) = \{\qiint{n}:n\in\NN_{0}\},$$ and every eigenvalue has infinite multiplicity. 
\end{enumerate}
\end{thmA}

As $q\to1^{-}$, $\Lam$ becomes $A^*A=(\Mzb-\pp)\ppb$, whose eigenspaces are the true polyanalytic Fock spaces $\Fock^{(n)}$; in this sense the operator $\Lam$ represents the $q$-deformed magnetic Laplacian. $\Lam$ is essentially self-adjoint on $\Poly$ see Corollary~\ref{cor:esa}. 

\begin{thmB}
For $m,n\in\NN_{0}$ set $H_{m,n}:=q^{\,mn}(A^*)^{n}z^{m}$, with the convention $H_{m,-1}=H_{-1,n}=0$.
\begin{enumerate}[label=\textup{(\roman*)},leftmargin=*]
\item $(A^*)^{n} =\bigl(\prod_{j=0}^{n-1}(\Mzb-q^{\,j}\Rq)\bigr)q^{-n\Nz}$ and
\begin{equation}
\label{eq:Hexpl}
  H_{m,n}=\sum_{k=0}^{m\wedge n}(-1)^{k}q^{\binom{k}{2}} \qbin{m}{k}\qbin{n}{k}\qfac{k}\;z^{\,m-k}\zb^{\,n-k}.
\end{equation}
\item $\Rqb H_{m,n}=\qint{n}H_{m,n-1}$, $\Rq H_{m,n}=\qint{m}H_{m-1,n}$, $a^*H_{m,n}=q^{-\frac{m+n}{2}}H_{m,n+1}$ and $\Lam H_{m,n}=\qiint{n}H_{m,n}$.
\item $\langle H_{m,n},H_{m',n'}\rangle
 =\delta_{m,m'}\delta_{n,n'}\,q^{\,mn}\qfac{m}\qfac{n}$ also
$\operatorname{span}\{H_{m,n}\}_{m,n\ge0}=\Poly$ and the $H_{m,n}$ form an orthogonal basis of the vector space $\Poly$.
\item For $|uz|<\frac{1}{1-q}$ and $|v\zb|<\frac{1}{1-q}$,
\begin{equation}
\label{eq:Hgen}
  \sum_{m,n\ge0}\frac{u^{m}v^{n}}{\qfac{m}\qfac{n}}\,H_{m,n}
  =\frac{\Eq(uz)\,\Eq(v\zb)}{\Eq(uv)} .
\end{equation}
\end{enumerate}
\end{thmB}

We relate our findings to the $q$-Hermite polynomials introduced by Ismail and Zhang~\cite{IsmailZhang}.

This allows us now to characterise our $q$-polyanalytic Fock spaces $\HHk{n}$.

\begin{thmC}
Put $h_{m,n}=H_{m,n}/\sqrt{q^{\,mn}\qfac{m}\qfac{n}}$ and let $\HHk{n}$ be the closure of $\operatorname{span}\{h_{m,n}\}_{m\ge0}$ in $\Lq$, so that $\HHk{0}=\HH$.
\begin{enumerate}[label=\textup{(\roman*)},leftmargin=*]
\item $\{h_{m,n}\}_{m\ge0}$ is an orthonormal basis of $\HHk{n}$, the spaces $\HHk{n}$ are mutually orthogonal and $\bigoplus_{n\ge0}\HHk{n}$ is the closure of $\Poly$ in $\Lq$.
\item $\sum_{m}|h_{m,n}(z)|^{2}$ converges for $|z|<r_{n}$ and diverges for $|z|\ge r_{n}$ off the finite sets $|z|^{2}\in\{s_{0},\dots,s_{n-1}\}$. Hence identifying $f=\sum_{m}c_{m}h_{m,n}\in\HHk{n}$ with the locally uniform limit of $\sum_{m}c_{m}h_{m,n}$ on $\DD_{r_n}$ (an injective identification, which agrees with $f$ almost everywhere on each circle $\TT_{r_j}$, $j>n$), $\HHk{n}$ is a reproducing kernel Hilbert space on $\DD_{r_{n}}$ with kernel $\sum_{m}h_{m,n}(z)\overline{h_{m,n}(w)}$ and on no larger disk and its elements are polyanalytic of order $n+1$ there. 
\item Let $\Bq$ be the $q$-Bargmann transform of Definition~\ref{def:qBarg}. Then $\|(a^*)^{n}\Bq\psi_{m}\|^{2}=q^{-\binom{n}{2}}\qfac{n}$ for every $m$, and
\begin{equation}\label{eq:transform}
  \mathcal{B}_{q}^{(n)}
  := \frac{1}{\sqrt{\qiint{n}!}}\bigl(a^*\bigr)^{n}\Bq= 
   \frac{q^{\frac{1}{2}\binom{n}{2}}}{\sqrt{\qfac{n}}}
   \bigl(a^*\bigr)^{n}\Bq ,
\end{equation}
defined on $\operatorname{span}\{\psi_{m}\}_{m\ge0}$, extends to a unitary map of $L^{2}(\RR)$ onto $\HHk{n}$ with $\mathcal{B}_{q}^{(n)}\psi_{m}=h_{m,n}$.
\end{enumerate}
\end{thmC}

Theorems A, B and~C are proved in Sections~\ref{sec:adjoint}--\ref{sec:levels}.

\section{Proof of Theorem A}\label{sec:adjoint}

\begin{proof}[Proof of Theorem \textup{A(i)}]
Set $A^*=(\Mzb-\Rq)q^{-\Nz}$, so that $A^*e_{c,d}=q^{-c}\bigl(e_{c,d+1}-\qint{c}e_{c-1,d}\bigr)$. In \eqref{eq:gram} the two resulting terms carry the same constraint, since $a-b=c-(d+1)$ and $a-b=(c-1)-d$ are both $a-b=c-d-1$. Hence
\begin{align*}
  \langle e_{a,b},A^*e_{c,d}\rangle
  &=q^{-c}\,\delta_{a-b,\,c-d-1}\,\qfac{a+d}
   \bigl(\qint{a+d+1}-\qint{c}\bigr)\\
  &=\delta_{a-b,\,c-d-1}\,\qint{b}\,\qfac{a+d}
  =\qint{b}\,\delta_{a-b+1,\,c-d}\,\qfac{a+d}
  =\langle Ae_{a,b},e_{c,d}\rangle ,
\end{align*}
where $\qint{a+d+1}-\qint{c}=\frac{q^{c}-q^{\,a+d+1}}{1-q} =q^{c}\qint{a+d+1-c}$ and $a+d+1-c=b$ on the constraint set. The monomials span $\Poly$ and adjoints are unique. For the commutator,
\begin{align}
  AA^*e_{a,b}&=q^{-a}\bigl(\qint{b+1}e_{a,b}
     -\qint{a}\qint{b}e_{a-1,b-1}\bigr),\label{eq:AAdag}\\
  A^*Ae_{a,b}&=q^{-a}\qint{b}\bigl(e_{a,b}
     -\qint{a}e_{a-1,b-1}\bigr),\label{eq:AdagA}
\end{align}
and
$ \left(AA^*-A^*A\right)e_{a,b}=q^{-a}\bigl(\qint{b+1}-\qint{b}\bigr)e_{a,b}=q^{\,b-a}e_{a,b}$.
\end{proof}

\begin{proof}[Proof of Theorem \textup{A(ii)}]
Write $J=b-a$ for the grading of $e_{a,b}$. Since $Ae_{a,b}$ lies in $\Pgr_{J-1}$ and both terms of $A^*e_{a,b}$ lie in $\Pgr_{J+1}$,
\[
  \widetilde a\,e_{a,b}=d(J-1)\qint{b}e_{a,b-1},
  \qquad
  \widetilde a^*e_{a,b}
  =\overline{d(J)}\,q^{-a}\bigl(e_{a,b+1}-\qint{a}e_{a-1,b}\bigr).
\]
Composing and using \eqref{eq:AAdag} and \eqref{eq:AdagA},
\begin{align}
  \widetilde a\widetilde a^*e_{a,b}
   &=|d(J)|^{2}q^{-a}\bigl(\qint{b+1}e_{a,b}
      -\qint{a}\qint{b}e_{a-1,b-1}\bigr),\label{eq:tildeaa}\\
  \widetilde a^*\widetilde a\,e_{a,b}
   &=|d(J-1)|^{2}\qint{b}\,q^{-a}\bigl(e_{a,b}
      -\qint{a}e_{a-1,b-1}\bigr).\label{eq:tildeaa2}
\end{align}
Equating the coefficients in the expression $\widetilde a\widetilde a^*-\lambda\widetilde a^*\widetilde a\,$ yields the following. First, the coefficient of $e_{a-1,b-1}$ in $\widetilde a\widetilde a^*-\lambda\,\widetilde a^*\widetilde a$ is
\[
-\qint{a}\qint{b}\,q^{-a}\bigl(|d(J)|^{2}-\lambda|d(J-1)|^{2}\bigr).
\]
 By Lemma~\ref{lem:inj}, $e_{a-1,b-1}$ and $e_{a,b}$ are distinct elements of the monomial basis of $\Pgr_{J}$ due to which coefficient vanishes, hence
\[
|d(J)|^{2}=\lambda\,|d(J-1)|^{2}\qquad\text{for all }J\in\ZZ.
\]
Second, by this relation the coefficient of $e_{a,b}$ is
\[
 q^{-a}\bigl(|d(J)|^{2}\qint{b+1}-\lambda|d(J-1)|^{2}\qint{b}\bigr)={} |d(J)|^{2}q^{-a}\bigl(\qint{b+1}-\qint{b}\bigr)=|d(J)|^{2}q^{J}.
\]
Consequently, $|d(J)|^{2}=q^{-J}$ and therefore
\[
\lambda=\frac{|d(J)|^{2}}{|d(J-1)|^{2}}=q^{-1}.
\]

Conversely let $d(J)=q^{-\frac{J}{2}}$. Then $a$ and \(a^*\) are mutually adjoint on $\Poly$, since $q^{-\frac{\J}{2}}$ is symmetric on  this space. Moreover, by $|d(J)|^{2}=q^{-J}$ for all $J$,  \eqref{eq:tildeaa} and \eqref{eq:tildeaa2} take the following form
\begin{align}
  a\,a^*e_{a,b}&=q^{-b}\bigl(\qint{b+1}e_{a,b}
    -\qint{a}\qint{b}e_{a-1,b-1}\bigr),\label{eq:aadag}\\
  a^*a\,e_{a,b}&=\qint{b}q^{\,1-b}\bigl(e_{a,b}
    -\qint{a}e_{a-1,b-1}\bigr).\label{eq:adaga}
\end{align}
For \(\lambda=q^{-1}\) the remaining coefficient simplifies to
\[
q^{-b}\bigl(\qint{b+1}-\qint{b}\bigr)=1,
\]
which yields \eqref{eq:AC} and completes the proof.
\end{proof}
\begin{proof}[Proof of Theorem \textup{A(iii)}]
Symmetry follows from $\langle a^*au,v\rangle=\langle au,av\rangle=\langle u,a^*av\rangle$, both equalities by Theorem~A(ii); in particular $\langle\Lam u,u\rangle=\|au\|^{2}\ge0$. Moreover, by using \eqref{eq:adaga} with $q^{\,1-b}\qint{b}=\qiint{b}$ we obtain the stated action of $\Lam$. Since both monomials on the right-hand side belong to $\Pgr_{b-a}$, the operator $\Lam$ preserves the grading. Consequently, if $\Lam w=\mu w$ with $w=\sum_{J}w_{J}$, $w_{J}\in\Pgr_{J}$, then $\Lam w_{J}=\mu w_{J}$ for every $J$, so it suffices to consider eigenvectors in a single sector. 

Fix $J\in\ZZ$ and put $b_{0}:=\max(0,J)$; the monomials $v_{b}:=e_{b-J,b}$, $b\ge b_{0}$, form a basis of $\Pgr_{J}$. Then $\Lam v_{b}=\qiint{b}v_{b}-\qiint{b}\qint{b-J}v_{b-1}$, which is upper bidiagonal in the ordered basis $(v_{b_{0}},v_{b_{0}+1},\dots)$ with diagonal entries $\qiint{b}$, $b\ge b_0$.  The superdiagonal coefficient $\qiint{b}\qint{b-J}$ vanishes at $b=b_{0}$. By \eqref{eq:reciprocal} the sequence $\bigl(\qiint{b}\bigr)_{b\ge b_{0}}$ are strictly increasing and pairwise distinct. 

Suppose that $\Lam w=\mu w$ for some non-zero $w\in\Pgr_{J}$.
We have $ w=\sum_{b=b_{0}}^{B}c_{b}v_{b}$ with $c_{B}\neq0$,
and comparing the coefficients of $v_{B}$, we obtain
\[
  \bigl(\qiint{B}-\mu\bigr)c_{B}=0.
\]
Hence $\mu=\qiint{B}$.

Conversely, let $n\ge b_{0}$ and seek an eigenvector
\[
  w_{n,J}=\sum_{b=b_{0}}^{n}c_{b}v_{b},
  \qquad c_{n}=1.
\]
The equation $\Lam w_{n,J}=\qiint{n}w_{n,J}$ is equivalent to
\[
  c_{b}
  =\frac{\qiint{b+1}\,\qint{b+1-J}}
         {\qiint{b}-\qiint{n}}\,c_{b+1},
  \qquad b=n-1,\ldots,b_{0},
\]
which determines $c_{n-1},\dots,c_{b_{0}}$ uniquely. Thus $\qiint{n}$ is an eigenvalue in every $\Pgr_{J}$ with $J\le n$ and its multiplicity is infinite.
\end{proof}

\begin{proposition}
\label{prop:braid}
On $\Poly$ one has $[A,\Rq]=0$, $[\Lam,\Rq]=0$ and $a\Rq=q^{-\frac{1}{2}}\Rq a$.
\end{proposition}

\begin{proof}  
The first identity is trivial: by \eqref{eq:jackson2var}, $\Rq A e_{a,b}=\qint{a}\qint{b}e_{a-1,b-1}=A\Rq e_{a,b}$. For the second, Theorem A(iii) gives
\[
  \Lam\Rq e_{a,b}=\qint{a}\qiint{b}  \bigl(e_{a-1,b}-\qint{a-1}e_{a-2,b-1}\bigr)=\Rq\Lam e_{a,b}.
\]
For the third, with $J=b-a$,
\[
  a\Rq e_{a,b}=\qint{a}\qint{b}q^{-\frac{J}{2}}e_{a-1,b-1},
  \qquad
  \Rq a\,e_{a,b}=\qint{a}\qint{b}q^{-\frac{J-1}{2}}e_{a-1,b-1},
\]
which completes the proof.
\end{proof}

As $q\to1^{-}$, the relation \eqref{eq:AC} becomes $[a,a^*]=I$ and Proposition~\ref{prop:braid} becomes $a\Rq=\Rq a$, which recovers the two commuting canonical pairs of the classical theory.

\section{Proof of Theorem B}\label{sec:hermite}

\begin{lemma}\label{lem:intertwine}
For every $j\in\ZZ$, on $\Poly$,
\[
  q^{-j\Nz}\Rq q^{\,j\Nz}=q^{\,j}\Rq\qquad\text{and}\qquad
  q^{-j\Nz}\Mzb q^{\,j\Nz}=\Mzb.
\]
\end{lemma}

\begin{proof}
$q^{-j\Nz}\Rq q^{\,j\Nz}e_{a,b} =q^{-j(a-1)}q^{\,ja}\qint{a}e_{a-1,b}=q^{\,j}\qint{a}e_{a-1,b}$ and $\Mzb$ does not change the $z$-degree.
\end{proof}

\begin{proof}[Proof of Theorem \textup{B(i)}]
Set $C=\Mzb-\Rq$ and $W=q^{-\Nz}$, so that $A^*=CW$. By Lemma~\ref{lem:intertwine} one has
\[
W^{\,j}CW^{-j}=\Mzb-q^{\,j}\Rq=:C_j, \quad \text{that is} \quad  W^{\,j}C=C_jW^{\,j}.
\]
Reordering factors gives the following
\[
(CW)^n = C\cdot\bigl(WCW^{-1}\bigr)\cdots \bigl(W^{\,n-1}CW^{-(n-1)}\bigr)W^n = C_0C_1\cdots C_{n-1}W^n.
\]
The $C_j$ commute with one another because $\Mzb$ and $\Rq$ act on independent variables. Moreover, since $q^{-n\Nz}z^m=q^{-mn}z^m$, it follows that
\[
H_{m,n}=\prod_{j=0}^{n-1}\bigl(\Mzb-q^{\,j}\Rq\bigr)z^m.
\]
By using the Gauss' formula \eqref{eq:gauss}, with the operators $x=\Mzb$ and $a=-\Rq$ , gives
\[
H_{m,n}=\sum_{k=0}^{n}(-1)^k q^{\binom{k}{2}}\qbin{n}{k}\Rq^{\,k}\Mzb^{\,n-k}z^m.
\] 
By using the identity $\Rq^{\,k}z^m=\frac{\qfac{m}}{\qfac{m-k}}z^{m-k}$ (which vanishes for $k>m$), followed by \eqref{eq:qbinfac}, we obtain \eqref{eq:Hexpl}. This completes the proof.
\end{proof}

Each term in the expansion \eqref{eq:Hexpl} belongs to $\Pgr_{n-m}$. The polynomial has leading term $z^{m}\zb^{\,n}$ and its coefficients are real-valued and symmetric with respect to $m$ and $n$. Moreover,
\begin{equation}\label{eq:Hconj}
  \overline{H_{m,n}}=H_{n,m},\qquad H_{m,0}=z^{m},\qquad
  H_{0,n}=\zb^{\,n}.
\end{equation}
For instance $H_{1,1}=z\zb-1$, $H_{1,2}=z\zb^{\,2}-(1+q)\zb$ and $H_{2,2}=z^{2}\zb^{\,2}-(1+q)^{2}z\zb+q(1+q)$.

\begin{proof}[Proof of Theorem \textup{B(ii)}]
Applying $\Rqb$ to \eqref{eq:Hexpl} gives;
\[
\Rqb H_{m,n} =  \sum_{k=0}^{m \wedge (n-1)}(-1)^k  q^{\binom{k}{2}} \qbin{m}{k} \qbin{n}{k} \qint{n-k} \qfac{k}\;z^{\,m-k}\zb^{\,(n-1)-k}
\]
Using
\[
  \qbin{n}{k}\,\qint{n-k}=\qint{n}\,\qbin{n-1}{k},
  \qquad 0\le k\le m\wedge(n-1),
\]
and \eqref{eq:Hexpl} with $n-1$ in place of $n$, we get $\Rqb H_{m,n}=\qint{n}H_{m,n-1}$.

Since $\overline{\Rq f}=\Rqb\,\overline{f}$ and $\overline{H_{m,n}}=H_{n,m}$ by \eqref{eq:Hconj} we get second identity.

For the third identity, the definition of $H_{m,n}$ gives 
\[
  A^*H_{m,n}=q^{\,mn}\bigl(A^*\bigr)^{n+1}z^{m}  =q^{\,mn}q^{-m(n+1)}H_{m,n+1}=q^{-m}H_{m,n+1},
\]
and since $H_{m,n}\in\Pgr_{n-m}$, $q^{-\frac{\J}{2}}$ acts on $H_{m,n}$ as multiplication by $q^{-\frac{n-m}{2}}$. Hence
\begin{equation}\label{eq:ladderup}
  a^*H_{m,n}=q^{-\frac{n-m}{2}}A^*H_{m,n}   =q^{-\frac{m+n}{2}}H_{m,n+1}.
\end{equation}
Finally, $AH_{m,n}=\qint{n}H_{m,n-1} \in \Pgr_{n-1-m}$, so
\begin{align*}
  \Lam H_{m,n}=a^*a\,H_{m,n}
  &=\qint{n}\,q^{-\frac{n-1-m}{2}}\,a^*H_{m,n-1}\\
  &=\qint{n}\,q^{-\frac{n-1-m}{2}}\,q^{-\frac{m+n-1}{2}}H_{m,n} \\
  &=\qint{n}\,q^{\,1-n}H_{m,n}\overset{\eqref{eq:reciprocal}}= \qiint{n} H_{m,n},
\end{align*}
which completes the proof.
\end{proof}

\begin{proof}[Proof of Theorem \textup{B(iii)}]
If \(n\neq n'\), then \(H_{m,n}\perp H_{m',n'}\), because these vectors are eigenvectors of the symmetric operator \(\Lam\) associated with distinct eigenvalues \(\qiint{n}\) and \(\qiint{n'}\). If \(n=n'\) and \(m\neq m'\), then \(H_{m,n}\in \Pgr_{n-m}\) and \(H_{m',n}\in \Pgr_{n-m'}\) are orthogonal by \eqref{eq:gram}.

 By the definition of $H_{m,n}$, we have
\[
H_{m,n}=q^{\,m}A^*H_{m,n-1}.
\]
Therefore, by Theorem A(i),
\[
  \|H_{m,n}\|^{2}
  =q^{\,m}\bigl\langle H_{m,n-1},AH_{m,n}\bigr\rangle
  =q^{\,m}\qint{n}\|H_{m,n-1}\|^{2}.
\]
Iterating this relation and by the initial condition \(\|H_{m,0}\|^{2}=\|z^{m}\|^{2}=\qfac{m}\) yields
\[
\|H_{m,n}\|^{2}=q^{\,mn}\qfac{m}\qfac{n},
\]
which is strictly positive. Moreover, since \(H_{m,n}\) is monic of bidegree \((m,n)\) and contains only monomials of strictly lower total degree within the same graded component, the resulting change-of-basis in each \(\Pgr_{J}\) is triangular with unit diagonal, hence it spans \(\Poly\), and, on each \(\Pgr_{J}\), it coincides with the monic Gram--Schmidt orthogonalisation of the monomial basis.
\end{proof}

\begin{proof}[Proof of Theorem \textup{B(iv)}]
Using \eqref{eq:Hexpl} together with \eqref{eq:qbinfac}, we obtain
\[
\frac{\qbin{m}{k}\qbin{n}{k}\qfac{k}}{\qfac{m}\qfac{n}} =\frac{1}{\qfac{k}\,\qfac{m-k}\,\qfac{n-k}}.
\]
By substituting  \(u^{m}=u^{k}u^{m-k}\) and \(v^{n}=v^{k}v^{n-k}\), \eqref{eq:Hgen} can be written as 
\[
\Bigl(\sum_{k\ge0}\frac{(-1)^{k}q^{\binom{k}{2}}(uv)^{k}}{\qfac{k}}
\Bigr)
\Bigl(\sum_{i\ge0}\frac{(uz)^{i}}{\qfac{i}}\Bigr)
\Bigl(\sum_{l\ge0}\frac{(v\zb)^{l}}{\qfac{l}}\Bigr)
\overset{\eqref{eq:qexp}}=
\frac{\Eq(uz)\,\Eq(v\zb)}{\Eq(uv)}.
\]
\end{proof}

\begin{remark}\label{rem:IZ}
Let $H^{\mathrm{IZ}}_{m,n}(z_{1},z_{2}\,|\,q)$ be the first family of two-dimensional $q$-Hermite polynomials of Ismail and Zhang~\cite[(3.1)]{IsmailZhang}. Comparing with \eqref{eq:Hexpl} and using $\qfac{k}=\rad^{2k}(q;q)_{k}$, one finds
\[
  H_{m,n}(z,\zb)=\rad^{\,m+n}\,H^{\mathrm{IZ}}_{m,n}\bigl(z/\rad,\zb/\rad\,\bigm|\,q\bigr),
\]
and the image of the orthogonality measure of~\cite[Theorem~3]{IsmailZhang} under $z\mapsto\rad z$, multiplied by $(q;q)_{\infty}$, is $\mu_{q}$. The present paper derives this family from the formal adjoint of $\Rqb$ in $\Lq$; this gives in addition the relation \eqref{eq:AC}, the operator $\Lam$ with its eigenvalue equation and the spaces $\HHk{n}$.

Our result in Theorem~C differs from~\cite{Arjika} as follows. In~\cite{Arjika} each level is treated separately: the space at level $n$ is the closure of $\operatorname{span}\{h_{j,n}\}_{j\ge0}$ in $L^{2}(\DD_{r_n},\mu_{q})$, a different Hilbert space for each $n$ and no operator having these spaces as eigenspaces is considered. Here all levels lie in $\Lq$; they are mutually orthogonal, their closed direct sum is $\Lq$ and they are the eigenspaces of $\Lam$ (Theorem~C(i) and Corollary~\ref{cor:esa}). The norm of $\HHk{n}$ is that of $\Lq$, not that of $L^{2}(\DD_{r_n},\mu_{q})$; for the latter the $h_{j,n}$ remain orthogonal but are not normalised and the two norms are not even equivalent. 
\end{remark}

\section{\texorpdfstring{Proof of Theorem C}{Proof of Theorems C} and the kernel of the \texorpdfstring{$q$}{q}-polyanalytic Fock space}\label{sec:levels}

\begin{proof}[Proof of Theorem \textup{C(i)}]
Part (i) of Theorem C is immediate from Theorem B(iii). The $h_{m,n}$ are orthonormal. The spaces $\HHk{n}$ are mutually orthogonal because the $H_{m,n}$ with distinct $n$ are orthogonal, and orthogonality passes to closed linear spans. The closed linear span of $\bigcup_{n}\HHk{n}$ contains $\operatorname{span}\{H_{m,n}\}_{m,n}=\Poly$ and is contained in $\overline{\Poly}$; hence it equals $\overline{\Poly}$, which is $\Lq$ by the density of $\Poly$.
\end{proof}

\begin{corollary}\label{cor:esa}
The operator $\Lam$ is essentially self-adjoint on $\Poly$ in $\Lq$. The spectrum of its closure is $\{\qiint{n}:n\in\NN_{0}\}$; each $\qiint{n}$ is an eigenvalue of infinite multiplicity, with eigenspace $\HHk{n}$.
\end{corollary}\begin{proof}
$\Lam$ is symmetric (Theorem~A(iii)) and densely defined (Section~\ref{sec:prelim}) and by Theorems~B(ii) and C(i) the family $\{h_{m,n}\}_{m,n\ge0}\subset\Poly$ is an orthonormal basis of $\Lq$ with $\Lam h_{m,n}=\qiint{n}h_{m,n}$. Let $v\in\Lq$ satisfy $\langle\Lam u,v\rangle=\langle u,\pm iv\rangle$ for all $u\in\Poly$. Taking $u=h_{m,n}$ gives $(\qiint{n}\pm i)\langle h_{m,n},v\rangle=0$ for all $m,n$, hence $v=0$. Thus $\Lam$ is essentially self-adjoint. The operator acting by $h_{m,n}\mapsto\qiint{n}h_{m,n}$ on its maximal domain is self-adjoint and extends $\Lam$, so it is the closure of $\Lam$. Since $\qiint{n}\nearrow\infty$ by \eqref{eq:reciprocal}, its spectrum is $\{\qiint{n}:n\in\NN_{0}\}$ and the eigenspace for $\qiint{n}$ is the closed span of $\{h_{m,n}\}_{m\ge0}$, that is, $\HHk{n}$.
\end{proof}

\begin{lemma}\label{lem:radialform}
Let $n\in\NN_{0}$ and $m\ge n$. Then $H_{m,n}=z^{\,m-n}P_{m,n}\bigl(|z|^{2}\bigr)$, where
\begin{equation}
\label{eq:radialform}
  P_{m,n}(s)=\sum_{k=0}^{n}(-1)^{k}q^{\binom{k}{2}}\qbin{m}{k}
   \qbin{n}{k}\qfac{k}\,s^{\,n-k}.
\end{equation}
Moreover $\qbin{m}{k}\le\frac{1}{(q;q)_{k}}$ for every $m$ and $P_{m,n}\to P_{\infty,n}$ as $m\to\infty$, where $P_{\infty,n}(s)=\prod_{j=0}^{n-1}(s-s_{j})$.
\end{lemma}

\begin{proof}
Assume that \(m\ge n\ge k\). Then one may factor
\[
z^{m-k}\zb^{\,n-k}=z^{\,m-n}(z\zb)^{\,n-k},
\]
and substituting this identity into \eqref{eq:Hexpl} implies \eqref{eq:radialform}. Furthermore, the \(q\)-binomial coefficient gives the following
\[
\qbin{m}{k}=\frac{(q^{\,m-k+1};q)_{k}}{(q;q)_{k}} =\frac{1}{(q;q)_{k}}\prod_{i=1}^{k}\bigl(1-q^{\,m-k+i}\bigr).
\]
When \(k\) is fixed in the limit \(m\to\infty\), the finite product \(\prod_{i=1}^{k}(1-q^{\,m-k+i})\) increases monotonically to \(1\). Hence \(\qbin{m}{k}\) increases to \(\frac{1}{(q;q)_{k}}\) and, in particular, is uniformly bounded above by this limiting value. In the case of \(m<k\), the same bound holds trivially, since \(\qbin{m}{k}=0\) by \eqref{eq:qbinfac}. Hence $\qbin{m}{k}\qfac{k}\to\frac{\qfac{k}}{(q;q)_{k}}=(1-q)^{-k}$, so the coefficient of $s^{\,n-k}$ in $P_{m,n}$ converges to $(-1)^{k}q^{\binom{k}{2}}\qbin{n}{k}(1-q)^{-k}$. By \eqref{eq:gauss} with $x=s$ and $a=-\frac{1}{1-q}$, 
\begin{align*}
    P_{\infty,n}(s) &=  \sum_{k=0}^{n}(-1)^{k}\qbin{n}{k}q^{\binom{k}{2}}\Bigl(\frac{1}{1-q}\Bigr)^{k}s^{\,n-k}  \\
     &\overset{\eqref{eq:gauss}} = \prod_{j=0}^{n-1}\Bigl(s-\frac{q^{\,j}}{1-q} \Bigr)
     = \prod_{j=0}^{n-1}\bigl(s- s_j\bigr).
\end{align*}
\end{proof}

\begin{proof}[Proof of Theorem \textup{C(ii)}]
Set \(s=|z|^{2}\) and assume $s>0$. For $m<n$ we have only finitely many terms with $m<n$. Hence, by Lemma~\ref{lem:radialform}, it is enough to study
\[
  |h_{m,n}(z,\bar z)|^2
  =\frac{|z|^{2(m-n)}|P_{m,n}(s)|^2}
  {q^{mn}\qfac{m}\qfac{n}}
  =\frac{s^{m-n}P_{m,n}(s)^2}
  {q^{mn}\qfac{m}\qfac{n}}
  =:T_m(s),\quad m\ge n .
\]
Using \(\qfac{m}=\frac{(q;q)_{m}}{(1-q)^{m}}\) together with the bounds \((q;q)_{\infty}\le(q;q)_{m}\le1\), we get constants \(0<c_{1}\le c_{2}<\infty\), depending only on \(q,n,s\), such that
\[
  c_{1}P_{m,n}(s)^{2}\theta^{\,m}\le T_{m}(s)
  \le c_{2}P_{m,n}(s)^{2}\theta^{\,m},
  \qquad \theta=\frac{s(1-q)}{q^{\,n}}=\frac{s}{s_{n}} .
\]
We analyse two conditions:
\begin{enumerate}
    \item If \(s<s_{n}\), then \(\theta<1\): Since \(s<s_{n}\le s_{j}\) for \(j\le n-1\), every factor in \(P_{\infty,n}(s)\) is negative and therefore we have \(P_{\infty,n}(s)\neq0\). As \(P_{m,n}(s)\to P_{\infty,n}(s)\), the  quantities \(P_{m,n}(s)^{2}\) remain bounded, implying \(\sum_{m}T_{m}(s)<\infty\).
     \item If \(s\ge s_{n}\) and \(s\notin\{s_{0},\dots,s_{n-1}\}\), then we have \(\theta\ge1\). Here, still we have \(P_{\infty,n}(s)\neq0\). It follows that \(T_{m}(s)\) is bounded from below by a positive constant for sufficiently large \(m\) and thus the series diverges. Since \(\{s_{0},\dots,s_{n-1}\}\) is finite, any disk properly containing \(\DD_{r_{n}}\) contains points where point evaluation is unbounded.
\end{enumerate}

Let $f=\sum_{m=0}^{\infty}c_mh_{m,n}\in\HHk{n}$, so that $\sum_{m=0}^{\infty}|c_m|^2=\|f\|^2$. Also, we can write $f$ by its Almansi components $F_\nu$, as
\[
      f(z)=\sum_{\nu=0}^{n}\bar z^\nu F_\nu(z),
  \qquad
  F_\nu(z):=\sum_{m=n-\nu}^{\infty}
  c_m\gamma_{m,\nu}z^{m-n+\nu}, \qquad 0\le\nu\le n,
\]
Using \eqref{eq:Hexpl}, we can write
\[
  h_{m,n}(z,\zb)
  =\sum_{\substack{0\le\nu\le n\\m\ge n-\nu}}
    \gamma_{m,\nu}z^{m-n+\nu}\zb^\nu
\]
where $\gamma_{m,\nu}$ is obtained from the term $k=n-\nu$ in \eqref{eq:Hexpl}. By using Lemma~\ref{lem:radialform} we have
\[
  |\gamma_{m,\nu}|
  \le
  \frac{\qbin{n}{n-\nu}\qfac{n-\nu}}
  {(q;q)_{n-\nu}\sqrt{q^{mn}\qfac{m}\qfac{n}}}
  =\frac{C_{n,\nu}}{\sqrt{q^{mn}\qfac{m}}},
\]
with $C_{n,\nu}$ independent of $m$. Since $|c_m|\le\|f\|$ and $\qfac{m}^{-1/m}\to1-q$, it follows that
\[
  \limsup_{m\to\infty}|c_m\gamma_{m,\nu}|^{1/m}  \le q^{-n/2}\sqrt{1-q}=r_n^{-1}.
\]
Hence $F_\nu(z)$ is analytic on $\DD_{r_n}$, while the fixed shift in the exponent does not affect the radius of convergence.
These series converge locally uniformly and therefore
\[
  \sum_{m=0}^{M}c_mh_{m,n}(z,\zb)
  \longrightarrow
  \sum_{\nu=0}^{n}\zb^\nu F_\nu(z)
\]
locally uniformly on $\DD_{r_n}$. By \eqref{eq:almansi}, this defines a polyanalytic representative of $f$ of order at most $n+1$.

Write $F:=\sum_{\nu=0}^{n}\zb^{\,\nu}F_{\nu}$ for this representative. The map $f\mapsto F$ is injective: by the uniqueness in \eqref{eq:almansi}, $F=0$ forces $F_{n}=0$ and the Taylor coefficients of $F_{n}$ are $c_{m}(q^{mn}\qfac{m}\qfac{n})^{-1/2}$, so all $c_{m}$ vanish. For $j>n$ the circle $\TT_{r_j}$ lies in $\DD_{r_n}$; there the partial sums converge uniformly to $F$ and, by the restriction bound of Section~\ref{sec:prelim}, in $L^{2}(\TT_{r_j},d\lambda_{r_j})$ to $f$, so $F=f$ almost everywhere on $\TT_{r_j}$. Finally, for $z\in\DD_{r_n}$,
\[
  F(z)=\Bigl\langle f,\sum_{m}\overline{h_{m,n}(z)}\,h_{m,n}\Bigr\rangle,
\]
the second argument being in $\HHk{n}$ by the convergence part; so point evaluation is bounded on $\DD_{r_n}$.
\end{proof}

For part (iii) of Theorem C we use a $q$-analogue of the Bargmann transform.

\begin{definition}\label{def:qBarg}
Let $\{\psi_{m}\}_{m\ge0}$ be the real-valued normalised Hermite functions. For $(z,x)\in\DD_{\rad}\times\RR$ let $\Aq(z,x)=\sum_{m\ge0}\frac{z^{m}}{\sqrt{\qfac{m}}}\psi_{m}(x)$ and
\[
  \Bq[\varphi](z)=\int_{\RR}\Aq(z,x)\varphi(x)\,dx .
\]

\end{definition}

\begin{theorem}\label{thm:qiso}
Let $0<q<1$ and $\rad=(1-q)^{-1/2}$. The series defining
$A_q(z,\cdot)$ converges in $L^2(\RR)$ if and only if $|z|<\rad$.
The convergence is uniform on compact subsets of $\DD_{\rad}$ in the
$L^2(\RR)$ norm, and
\[
  \|A_q(z,\cdot)\|_{L^2(\RR)}^2=\Eq(|z|^2)
  \qquad (|z|<\rad).
\]
For every $\varphi\in L^2(\RR)$ and $z\in\DD_{\rad}$,
\[
  \Bq[\varphi](z)
  =\int_{\RR} \Aq(z,x)\varphi(x)\,dx
  =\sum_{m=0}^{\infty}\langle\varphi,\psi_m\rangle e_{m,q}(z).
\]
With realization of $\HH$, $B_q:L^2(\RR)\to\HH$ is unitary and $B_q\psi_m=e_{m,q}$ for every
$m\ge0$.
\end{theorem}

\begin{proof}
By orthonormality the partial sums are Cauchy in $L^{2}(\RR)$ exactly when $\sum_{m}\frac{|z|^{2m}}{\qfac{m}}=\Eq(|z|^{2})<\infty$, that is for $|z|<\rad$. Since $(q;q)_{\infty}\le(q;q)_{m}\le1$, the terms $\frac{|z|^{2m}}{\qfac{m}}=\frac{((1-q)|z|^{2})^{m}}{(q;q)_{m}}$ decay geometrically if $|z|<\rad$ and are at least $1$ if $|z|\ge\rad$; for $|z|\le r<\rad$ the tails are dominated by those at $r$, which gives the uniformity. Cauchy--Schwarz then justifies term-by-term integration and
\[
  \Bq\psi_{m}=\sum_{l}\frac{z^{l}}{\sqrt{\qfac{l}}} \int_{\RR}\psi_{l}\psi_{m}=e_{m,q},
\]
so the map carries an orthonormal basis onto one. More precisely, since $\{\psi_{m}\}$ and $\{e_{m,q}\}$ are orthonormal bases of $L^{2}(\RR)$ and $\HH$, the map $\varphi\mapsto\sum_{m}\langle\varphi,\psi_{m}\rangle e_{m,q}$ is unitary from $L^{2}(\RR)$ onto $\HH$ and term-by-term integration shows that $\Bq\varphi=\sum_{m}\langle\varphi,\psi_{m}\rangle e_{m,q}$ on $\DD_{\rad}$, which is the realisation of this element of $\HH$ described in Section~\ref{sec:prelim}.
\end{proof}

As $q\to1^{-}$, the Hermite-series representative of $\Aq$
converges pointwise to the classical Bargmann kernel~\cite{Barg},
\[
  A_1(z,x)=\pi^{-1/4}
  \exp\!\left(-\tfrac12(z^2+x^2)+\sqrt{2}\,zx\right).
\]
For the classical limit of a related transform constructed from
continuous $q$-Hermite functions, see~\cite[Corollary~1.1]{Arjika}.

At $q=0$, the target space is the classical Hardy space $H^2(\DD)$
\cite[Section~6]{Arik}. Setting $[m]_0!=1$ in
Definition~\ref{def:qBarg} gives $B_0\psi_m=z^m$. Since the monomials
form an orthonormal basis of $H^2(\DD)$, the resulting map
$B_0\colon L^2(\RR)\to H^2(\DD)$ is unitary.

\begin{proof}[Proof of Theorem \textup{C(iii)}]
By Theorem B(ii), $a^*H_{m,l}=q^{-\frac{m+l}{2}}H_{m,l+1}$; iterating over $0\le l\le n-1$ and using $H_{m,0}=z^{m}$, we get
\[
  \bigl(a^*\bigr)^{n}z^{m}
  =\Bigl(\prod_{l=0}^{n-1}q^{-\frac{m+l}{2}}\Bigr)H_{m,n}
  =q^{-\frac{1}{2}\left(mn+\binom{n}{2}\right)}H_{m,n}.
\]
Dividing by $\sqrt{\qfac{m}}$ and using Theorem B(iii),
\[
  \bigl\|\bigl(a^*\bigr)^{n}e_{m,q}\bigr\|^{2}
  =q^{-\left(mn+\binom{n}{2}\right)}
   \frac{q^{\,mn}\qfac{m}\qfac{n}}{\qfac{m}}
  =q^{-\binom{n}{2}}\qfac{n}{}=\qiint{n}!,
\]
which is independent of $m$.

Since $\Bq\psi_{m}=e_{m,q}$ by Theorem~\ref{thm:qiso},
\[
  \mathcal{B}_{q}^{(n)}\psi_{m}
  =\frac{q^{\frac{1}{2}\binom{n}{2}}}{\sqrt{\qfac{n}}}\,
   q^{-\frac{1}{2}\left(mn+\binom{n}{2}\right)}
   \frac{H_{m,n}}{\sqrt{\qfac{m}}}
  =\frac{H_{m,n}}{\sqrt{q^{\,mn}\qfac{m}\qfac{n}}}=h_{m,n},
\]
so $\mathcal{B}_{q}^{(n)}$ carries the orthonormal basis $\{\psi_{m}\}$ of $L^{2}(\RR)$ onto the orthonormal basis $\{h_{m,n}\}$ of $\HHk{n}$ and therefore extends from $\operatorname{span}\{\psi_{m}\}_{m\ge0}$ to a unitary.
\end{proof}

Finally, the reproducing kernel of $\HHk{1}$, defined on $\DD_{r_{1}}\times\DD_{r_{1}}$ by Theorem C(ii).

\begin{theorem}\label{thm:K1}
For $|z\wb|<s_{1}=\frac{q}{1-q}$,
\begin{equation}\label{eq:K1}
  K_{q}^{(1)}(z,w)=\sum_{m\ge0}h_{m,1}(z)\overline{h_{m,1}(w)}
  =\frac{1}{q}\,\Eq\Bigl(\frac{z\wb}{q}\Bigr)
   \bigl(1+z\wb+q\zb w-|z|^{2}-|w|^{2}\bigr).
\end{equation}
As $q\to1^{-}$ the right-hand side tends to $e^{z\wb}\bigl(1-|z-w|^{2}\bigr)$, the kernel of $\Fock^{(1)}$ (see, ~\cite{AF, Mouayn}).
\end{theorem}

\begin{proof}
    By using Theorem~B(iii) with polynomial $H_{m,n}$ defined in \eqref{eq:Hexpl}, we have
\[
  H_{m,1}(z,\zb)=z^m\zb-\qint{m}z^{m-1},
  \qquad
  \|H_{m,1}\|^2=q^m\qfac{m}.
\]
Thus, by substituting $u=z\wb$ we obtain
\[
\begin{aligned}
  K_q^{(1)}(z,w)
  &=
  \sum_{m=0}^{\infty}
  \frac{H_{m,1}(z,\zb)\overline{H_{m,1}(w,\wb)}}
       {q^m\qfac{m}}\\
       &= \zb w\underbrace{
    \sum_{m=0}^{\infty}\frac{u^m}{q^m\qfac{m}}
  }_{=\Eq\left(\frac{u}{q}\right)}
  -\bigl(|z|^2+|w|^2\bigr)\underbrace{
    \sum_{m=1}^{\infty}
    \frac{\qint{m}u^{m-1}}{q^m\qfac{m}}
  }_{=\frac1q\Eq\left(\frac{u}{q}\right)}
  +\underbrace{
    \sum_{m=1}^{\infty}
    \frac{\qint{m}^2u^{m-1}}{q^m\qfac{m}}
  }_{=\frac{1+u}{q}\Eq\left(\frac{u}{q}\right)}.
\end{aligned}
\]
All three series converge absolutely for $|u|<s_{1}$. Combining these identities, we conclude that
\[
  K_q^{(1)}(z,w)
  =\frac{1+z\wb+q\zb w-|z|^2-|w|^2}{q}
   \Eq\Bigl(\frac{z\wb}{q}\Bigr).
\]
For fixed $z,w$, letting $q\to1^{-}$ and using $\Eq(z\wb/q)\to e^{z\wb}$ (see, e.g.,~\cite{KC}) and $z\wb+\zb w-|z|^{2}-|w|^{2}=-|z-w|^{2}$ gives the limit.
\end{proof}

\section*{Acknowledgements}

 This study received financial support from a grant issued by the Methusalem programme under the Special Research Fund (BOF) at Ghent University (Grant number 01M01021). The third author is also supported by EPSRC grant UKRI3645. The second author was supported by CIDMA under the Portuguese Foundation for Science and Technology, Grants UID/04106/2025 and UID/PRR/04106/2025. Additionally, the first author's work is supported by the Higher Education Commission of Pakistan through the Overseas Scholarship Program (Phase III).

\section*{Use of generative AI}
We used GenAI (ChatGPT and Claude) to identify possible gaps in selected proofs and suggest intermediate steps. We also used it to compare our results with earlier work and to edit the exposition. We take responsibility for the proofs, the accuracy of the references and the final text.


\begin{thebibliography}{99}
\bibitem{AbreuACHA}
L.~D. Abreu,
\emph{Sampling and interpolation in Bargmann--Fock spaces of
polyanalytic functions}, Appl. Comput. Harmon. Anal. \textbf{29} (2010), no.~3, 287--302; MR2672228.

\bibitem{AF}
L.~D. Abreu and H.~G. Feichtinger,
\emph{Function spaces of polyanalytic functions},
in: Harmonic and Complex Analysis and its Applications, 1--38, Trends Math., Birkh\"auser/Springer, Cham, 2014; MR3203099.

\bibitem{ACKS}
D.~Alpay, P.~Cerejeiras, U.~K\"{a}hler and B.~Schneider,
\emph{Generalized $q$-Fock spaces and structural identities},
Proc. Amer. Math. Soc. \textbf{152} (2024), no.~6, 2463--2478; MR4741242.

\bibitem{AlpayPorat}
D.~Alpay and M.~Porat,
\emph{Generalized Fock spaces and the Stirling numbers},
J. Math. Phys. \textbf{59} (2018), no.~6, 063509, 12 pp.; MR3816055.

\bibitem{Arik}
M.~Arik and D.~D.~Coon,
\emph{Hilbert spaces of analytic functions and generalized coherent
states}, J. Math. Phys. \textbf{17} (1976), no.~4, 524--527; MR0495971.

\bibitem{Arjika}
S.~Arjika, O.~El Moize and Z.~Mouayn,
\emph{Une $q$-d\'{e}formation de la transformation de Bargmann
vraie-polyanalytique},
C. R. Acad. Sci. Paris S\'{e}r. I \textbf{356} (2018), no.~8, 903--910; MR3851546.
\bibitem{Mouayn}
N.~E.~Askour, A.~Intissar and Z.~Mouayn,
\emph{Explicit formulas for reproducing kernels of generalized Bargmann
spaces of $\CC^{n}$},
J. Math. Phys. \textbf{41} (2000), no.~5, 3057--3067,
doi:10.1063/1.533312.
\bibitem{Barg}
V.~Bargmann,
\emph{On a Hilbert space of analytic functions and an associated
integral transform, Part I},
Comm. Pure Appl. Math. \textbf{14} (1961), 187--214; MR0157250.

\bibitem{ElMoizeMouayn}
O.~El Moize and Z.~Mouayn,
\emph{A $q$-deformation of true-polyanalytic Bargmann transforms when
$q^{-1}>1$}, C. R. Math. Acad. Sci. Paris \textbf{359} (2021), no.~10, 1295--1305.

\bibitem{GR}
G.~Gasper and M.~Rahman,
\emph{Basic Hypergeometric Series}, 2nd ed.,
Encyclopedia of Mathematics and its Applications \textbf{96},
Cambridge University Press, 2004; MR2128719.

\bibitem{IsmailZhang}
M.~E.~H. Ismail and R.~Zhang,
\emph{On some $2D$ orthogonal $q$-polynomials},
Trans. Amer. Math. Soc. \textbf{369} (2017), no.~10, 6779--6821; MR3683094.

\bibitem{Ito}
K.~It\^{o},
\emph{Complex multiple Wiener integral},
Jap. J. Math. \textbf{22} (1952), 63--86; MR0063609.

\bibitem{jackson}
F.~H. Jackson,
\emph{On $q$-Functions and a Certain Difference Operator},
Trans. Roy. Soc. Edinburgh \textbf{46} (1909), 253--281 (pre-MR).

\bibitem{JungSato}
Y.-G.~Jung and R.~Sato,
\emph{$q$-deformed polyanalytic Ginibre point processes: construction and central limit theorems for linear statistics},
preprint, arXiv:2609.06904 (2026).

\bibitem{KC}
V.~Kac and P.~Cheung,
\emph{Quantum Calculus},
Universitext, Springer, New York, 2002; MR1865777.

\bibitem{LM}
H.~van Leeuwen and H.~Maassen,
\emph{A $q$ deformation of the Gauss distribution},
J. Math. Phys. \textbf{36} (1995), no.~9, 4743--4756,
doi:10.1063/1.530917.

\bibitem{ElMoizeMouaynRS}
Z.~Mouayn and O.~El Moize,
\emph{A set of $q$-coherent states for the Rogers--Szeg\H{o} oscillator},
Lett. Math. Phys. \textbf{111} (2021), Paper No.~143; arXiv:2102.10714.

\bibitem{RT16}
M.~Ruzhansky and N.~Tokmagambetov,
\emph{Nonharmonic analysis of boundary value problems},
Int. Math. Res. Not. IMRN \textbf{2016} (2016), no.~12, 3548--3615.

\bibitem{RT17a}
M.~Ruzhansky and N.~Tokmagambetov,
\emph{Wave equation for operators with discrete spectrum and irregular propagation speed},
Arch. Ration. Mech. Anal. \textbf{226} (2017), no.~3, 1161--1207.

\bibitem{RT17b}
M.~Ruzhansky and N.~Tokmagambetov,
\emph{Very weak solutions of wave equation for Landau Hamiltonian with irregular electromagnetic field},
Lett. Math. Phys. \textbf{107} (2017), no.~4, 591--618.

\bibitem{RT19a}
M.~Ruzhansky and N.~Tokmagambetov,
\emph{Wave equation for 2D Landau Hamiltonian},
Appl. Comput. Math. \textbf{18} (2019), no.~1, 69--78.

\bibitem{RT19b}
M.~Ruzhansky and N.~Tokmagambetov,
\emph{On nonlinear damped wave equations for positive operators. I. Discrete spectrum},
Differential Integral Equations \textbf{32} (2019), no.~7-8, 455--478; MR3945764.


\bibitem{Vas}
N.~L. Vasilevski,
\emph{Poly-Fock spaces},
in \emph{Differential Operators and Related Topics (Odessa, 1997)},
Vol.~I, Oper. Theory Adv. Appl. \textbf{117}, Birkh\"auser, Basel,
2000, pp.~371--386; MR1764974.


\end{thebibliography}
\end{document}